\documentclass[10pt,a4paper]{article}
\usepackage{amssymb,cite}
\usepackage{enumerate,verbatim}
\usepackage{amsmath}
\usepackage{amsfonts,mathrsfs}
\usepackage{amsthm,wasysym,tikz}
\usepackage{epsfig,lpic}
\usepackage[applemac]{inputenc}
\usepackage[english]{babel}

\usepackage{pdfsync}
\usepackage{xcolor,hyperref}
\usepackage[T1]{fontenc}
\usepackage{mathtools}

\usepackage[width=1\textwidth]{caption}

\definecolor{green}{rgb}{0,0.5,0}

\hypersetup{backref,colorlinks=true,linkcolor=blue,citecolor=blue}

\usepackage{marvosym}

\newcounter{teoremaganso}
\newcounter{appendix}

\newcounter{coryganso}

\renewenvironment{abstract}{\small\quotation\noindent
 {\bfseries \abstractname .}}{\endquotation \par}

\newenvironment{trivlistparm}[2]{\trivlist
\item[{#1 #2}]}{\endtrivlist}

\newenvironment{prooftext}[1]{\trivlistparm{\bfseries}{#1}}{\Qed\endtrivlistparm}

\catcode`\@=11

\def\resetthefootnote{\renewcommand{\thefootnote}{\@arabic\c@footnote} }
\def\@principiremex#1{\trivlist
 \item[\hskip \labelsep{\bfseries #1\ \thetheo.}]\ignorespaces}
\def\opar@principiremex#1[#2]{\trivlist
 \item[\hskip \labelsep{\bfseries #1\ \thetheo\ (#2).}]\ignorespaces}

\newcommand{\newTHEOremrom}[2]{\newenvironment{#1}{\refstepcounter{theo}\@ifnextchar[{\opar@principiremex{#2}}
{\@principiremex{#2}}}{\qedB\endtrivlist}} \catcode`\@=12
\DeclareMathSymbol{\square}{\mathord}{AMSa}{"03}
\newcommand{\qedB}{\nopagebreak\hspace*{\fill}$\square$\par}
\newcommand{\Qed}{\nopagebreak\hspace*{\fill}{\vrule width6pt height6pt depth0pt}\par}

\newtheorem {bigtheo} [teoremaganso] {Theorem}

\newTHEOremrom {defi} {Definition}
\newTHEOremrom {obs} {Remark}
\newTHEOremrom {ex} {Example}

\newcommand{\refc}[1]{\mbox{$(\ref{#1})$}}

\newcommand{\teoc}[1]{Theorem~\ref{#1}}

\newcommand{\figc}[1]{Figure~\ref{#1}}

\newcommand{\N}{\ensuremath{\mathbb{N}}}

\newcommand{\R}{\ensuremath{\mathbb{R}}}

\newcommand{\T}{\mathcal{T}}

\newcommand{\RP}{\ensuremath{\mathbb{RP}}}

\newcommand{\blue}[1]{{\color{blue}#1}}

\title{\bf Non-bifurcation of critical periods from \\ semi-hyperbolic polycycles of quadratic centers 
\footnotetext{This work has been partially funded by the Ministry of Science, Innovation and Universities of Spain through the grants PGC2018-095998-B-I00 and MTM2017-86795-C3-2-P, and by the Agency for Management of University and Research Grants of Catalonia through the grants 2017SGR1725 and 2017SGR1617.
}    
}
\author{D. Mar\'{\i}n,  M. Saavedra and J. Villadelprat
\\*[.1truecm]
{\small \textsl{Departament de Matem{\`a}tiques, Edifici Cc,
Universitat Aut{\`o}noma de Barcelona,}}\\*[-.05truecm]
{\small\textsl{08193 Cerdanyola del Vall\`es (Barcelona), Spain}}
\\*[-.05truecm]
{\small \textsl{Centre de Recerca Matem\`atica, Edifici Cc, Campus de Bellaterra,}}\\*[-.05truecm]
{\small \textsl{08193 Cerdanyola del Vall\`es (Barcelona), Spain}}
\\*[.1truecm]    
{\small \textsl{Departamento de Matem\'atica, Facultad de Ciencias F\'isicas y Matem\'aticas}}
\\*[-.05truecm]
{\small \textsl{Universidad de Concepci\'on, Barrio Universitario, Concepci\'on, Casilla 160-C, Chile}}
\\*[.1truecm]
{\small \textsl{Departament d'Enginyeria Inform\`{a}tica i Matem\`{a}tiques,
                        ETSE,}}
\\*[-.05truecm]
{\small \textsl{Universitat Rovira i Virgili, 43007
                     Tarragona, Spain}}}

\begin{document}
\maketitle
\abstract{
In this paper we consider the unfolding of saddle-node
\[
X= \frac{1}{xU_a(x,y)}\Big(x(x^\mu-\varepsilon)\partial_x-V_a(x)y\partial_y\Big),
\]
parametrized by $(\varepsilon,a)$ with $\varepsilon\approx 0$ and $a$ in an open subset $A$ of $\R^{\alpha},$
and we study the Dulac time $\T(s;\varepsilon,a)$ of one of its hyperbolic sectors. We prove (\teoc{A}) that the derivative $\partial_s\T(s;\varepsilon,a)$ tends to $-\infty$ as $(s,\varepsilon)\to (0^+,0)$ uniformly on compact subsets of $A.$ This result is addressed to study the bifurcation of critical periods in the Loud's family of quadratic centers. In this regard we show (\teoc{regular}) that no bifurcation occurs from certain semi-hyperbolic polycycles.}

\section{Introduction and main results}

The present paper deals with planar polynomial ordinary differential systems and we study the qualitative properties of the period function of centers. A singular point of a planar differential system is a \emph{center} if it has a punctured neighbourhood that consists entirely of periodic orbits surrounding it. The largest neighbourhood with this property is called the \emph{period annulus} of the center and we denote it by $\mathscr P.$ The \emph{period function} assigns to each periodic orbit in $\mathscr P$ its period. If the period function is constant then the center is called \emph{isochronous}. The study of the period function is a nontrivial problem and questions related to its behaviour have been extensively studied. Let us quote, for instance, the problems of isochronicity (see \cite{Colin,CMV,Xingwu}), monotonicity (see \cite{Chicone1,V07,Zhao}) or bifurcation of critical periods (see \cite{Chicone2,Gasull,GV10}). Aside from the intrinsic interest of these problems, the study of the period function is also important in the analysis of nonlinear boundary value problems and in perturbation theory. Indeed, for instance, under the condition of non-criticality of the period function, zeros of appropriate Melnikov functions guarantee the persistence of subharmonic periodic orbits of a Hamiltonian system after a small periodic non-autonomous perturbation (see \cite{Fonda,GH}). Most of the work on planar polynomial differential systems, including the present paper, is related to questions surrounding the well known Hilbert’s 16th problem (see \cite{Pugh,GMM,Roussarie,Ye} and references therein) and its various weakened versions.

Chicone~\cite{Chicone3} has conjectured that if a quadratic differential
system has a center with a period function which is not monotonic then, by an affine
transformation and a constant rescaling of time, it can be brought to the Loud normal
form
 \begin{equation}\label{loudB}
  \left\{
   \begin{array}{l}
    \dot u=-v+Buv, \\ \dot v=u+Du^2+Fv^2,
   \end{array}
  \right.
 \end{equation}
and that the period function of these centers has at most two critical periods. In fact,
there is much analytic evidence that the conjecture is true (see~\cite{CG,V07,Zhao} for
instance). On the other hand, it is proved in~\cite{GGV} that if $B=0$ then the period
function of the center at the origin of system~\refc{loudB} is globally monotonous. So,
from the point of view of the study of the period function, the most interesting stratum
of quadratic centers is the family~\refc{loudB} with $B\neq 0,$ which can be brought to
$B=1$ by means of a rescaling. Thus, using the vector field notation, in this paper we consider
\begin{equation}\label{loud0}
L_a\!:=(-v+uv)\partial_u+(u+Du^2+Fv^2)\partial_v\text{ where $a\!:=(D, F)\in \R^2$.}
\end{equation}
Following the terminology in~\cite{Chicone2}, we shall refer to this family as the \emph{dehomogenized Loud's centers}.

Compactifying $\R^2$ to the Poincaré disc, see for instance \cite{ADL}, the boundary of the period annulus~$\mathscr P$ of the center has two connected components, the center itself and a polycycle. We call them, respectively, the \emph{inner} and \emph{outer boundary} of the period annulus. Since period function is defined on the set of periodic orbits in~$\mathscr P,$  usually the first step is to parametrize this set, let us say $\{\gamma_s\}_{s\in
(0,1)},$ so that one can study the qualitative properties of the period function by
means of the map $s\mapsto\mbox{period of $\gamma_s$},$ which is analytic on~$(0,1).$
The \emph{critical periods} are the critical points of this function and its number,
character (maximum or minimum) and distribution do not depend on the particular
parametrization of the set of periodic orbits used. The dehomogenized Loud's family~\refc{loud0} depends on a two-dimensional parameter $a$ and our aim is to decompose $\R^2=\!\cup V_i$ so
that if~$a_1$ and~$a_2$ belong to the same set~$V_i,$ then the corresponding period
functions are qualitatively the same (i.e., their critical periods are
equal in number, character and distribution.) A parameter $a_0\in\R^2$ is a
\emph{regular value} if it belongs to the interior of some~$V_i,$ otherwise it is a
\emph{bifurcation value}. The set of bifurcation values is $\mathscr B\!:=\cup\,\partial V_i$
and, roughly speaking, it consists of those parameters $a_0\in\R^2$ for which some
critical period emerges or disappears as~$a\to a_0.$ There are three different situations to consider:
\begin{enumerate}[$(a)$]
\item Bifurcations of critical periods from the inner boundary (i.e., the center).
\item Bifurcations of critical periods from the interior of the period annulus.
\item Bifurcations of critical periods from the outer boundary (i.e., the polycycle).
\end{enumerate}
We refer the reader to \cite{MMV2} for the definition of these notions. 

\begin{figure}[t]
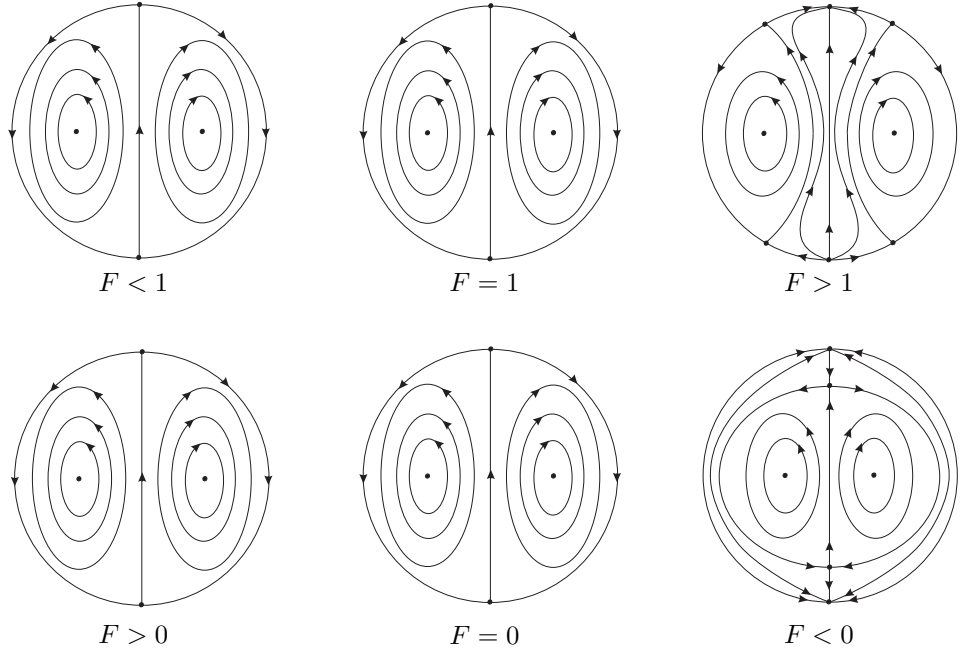

 \centering
 \begin{lpic}[l(0mm),r(0mm),t(0mm),b(5mm)]{fig1b(1.1)}
   \lbl[l]{53,45;$F=1$}
   \lbl[l]{11,45;$F<1$}
   \lbl[l]{93,45;$F>1$}
   \lbl[l]{53,2.5;$F=0$}
   \lbl[l]{11,2.5;$F>0$}
   \lbl[l]{93,2.5;$F<0$}    
 \end{lpic}
\caption{Phase portrait of $L_a$ in the Poincaré disc for $D\in(-1,0)$ and $F\approx 1$ (top) and for $D\in(-1,0)$ and $F\approx 0$ (bottom), where the center at the origin is shifted to the left for convenience and the vertical line is $u=1.$ The saddle-node singularity occurs for $F=1$ and $F=0$, respectively, and is placed at the intersection between the line at infinity (the boundary of the disc) and $u=1.$}\label{fig1}
\end{figure}

With regard to the dehomogenized Loud's centers \refc{loud0}, the bifurcation from the center was already solved by Chicone and Jacobs \cite{Chicone2}. Our goal is to study the bifurcation from the polycycle and to this end, together with P.~Marde\v si\'c, we have devoted a series of papers (see \cite{MMV,MMV2,MMV3,MMSV,MV,Saa20}). The polycycle consists of regular trajectories and singular points with a hyperbolic sector, which after the desingularization process give rise to hyperbolic saddles and saddle-nodes. Most of the cases studied so far correspond to hyperbolic polycycles, i.e., such that all the singularities at its vertices are hyperbolic saddles. Although this is the generic case in the family under consideration, in order to solve the problem we must tackle the non-hyperbolic polycycles as well. Among them there are two cases in which the polycycle has saddle-nodes, namely $(D,F)\in[-1,0]\times\{1\}$ and $(D,F)\in[-1,0]\times\{0\}$. In both cases the saddle-node bifurcation occurs at infinity, so one needs to extend the vector field to infinity by using, for instance, the Poincaré compactification (see \figc{fig1}). 
We treated the first case in \cite{MMSV}, where we proved a general result addressed to study the local passage 
through a singularity unfolding a saddle-node bifurcation. In the present paper we study the second case by adapting the general tools obtained in \cite{MMSV}. Let us briefly explain the similarities and differences between both cases. The key point is that the vector field $L_a$ has a Darboux first integral, which enables to find \emph{local} changes of coordinates that bring each unfolding to a suitable normal form. In the first case, see the proof of \cite[Theorem C]{MMSV}, the saddle-node unfolding $L_a$ for $D\in (0,1)$ and $F\approx 1$ can be brought locally to 
\[
\frac{1}{yU_a(x,y)}\left(x(x^2-\varepsilon)\frac{\partial}{\partial x}-V_a(x)y\frac{\partial}{\partial y}\right),\text{ where }\varepsilon=2(F-1)
\]
and $y=0$ corresponds to the line at infinity. (In its regard we remark that the polar factor can be neglected to draw the phase portrait but this cannot be done to study the time.) In this case the hyperbolic saddles at $\partial\mathscr P$ bifurcating from the saddle-node are placed at infinity for $F<1$ and $F>1$ (see the three phase portraits at the top of \figc{fig1}). In the present paper, by using the same techniques, we will show that the saddle-node unfolding $L_a$ for $D\in (0,1)$ and $F\approx 0$ can be brought locally to
\[
\frac{1}{xU_a( x, y)}\left(x( x^2-\varepsilon)\frac{\partial}{\partial  x}-V_a(x) y\frac{\partial}{\partial  y}\right),\text{ where }\varepsilon=-2F.
\]
Here the line at infinity corresponds to $x=0$. Unlike the previous case, the hyperbolic saddles at $\partial\mathscr P$ bifurcating from the saddle-node are located at infinity for $F>0$ but they are not for $F<0$ (see the three phase portraits at the bottom of \figc{fig1}). Thus, besides the  saddle-node bifurcation, in this case there is a second geometric phenomenon, namely, that the hyperbolic saddles at $\partial\mathscr P$ which bifurcate from the saddle-node located at infinity come to the finite plane for $F<0$. In this paper we deal with this more intrincate case and our main result states that there is no bifurcation of critical periods.

\begin{figure}[t]
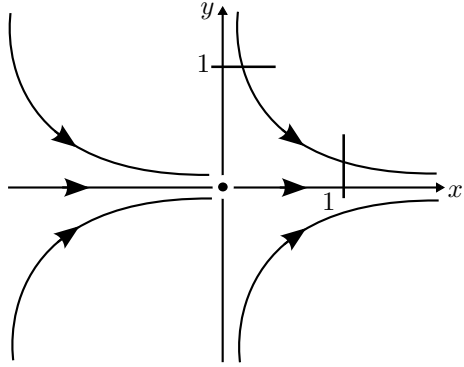

 \centering
 \begin{lpic}[l(0mm),r(0mm),t(0mm),b(5mm)]{temps_sella_node(1)}
   \lbl[l]{25,40;$1$}
   \lbl[l]{41.5,21.75;$1$} 
   \lbl[l]{58,23;$x$}        
   \lbl[l]{25.5,47;$y$}       
  \end{lpic}
 \caption{Transverse sections associated to the Dulac time $\T(s;\varepsilon,a)$ of the saddle-node unfolding \refc{X} for $\varepsilon=0$ (and taking $\mu$ odd).}\label{temps_SN}  
\end{figure}

In order to present our results in its full generality we adopt the framework  introduced in \cite{MMSV}.  We consider the unfolding of saddle-node  
\begin{equation}\label{X}
X= \frac{1}{xU_a(x,y)}\Big(x(x^\mu-\varepsilon)\partial_x-V_a(x)y\partial_y\Big),
\end{equation}
parametrized by $(\varepsilon,a)$ with $\varepsilon\approx 0$ and $a$ in an open subset $A$ of $\R^{\alpha}$, and where
\begin{itemize}
\item $\mu\in\N,$ 
\item $(x,y,a)\mapsto U_a(x,y)$ is analytic on $[-r,r]^2\times A$. Moreover, for each $a\in A,$ $U_a(0,0)>0$ and the Taylor series of $U_a(x,y)$ at $(0,0)$ is absolutely convergent on $[-r,r]^2.$ 
\item $(x,a)\mapsto V_{a}(x)$ is analytic on $[-r,r]\times A$ and, for all $a\in A,$ $V_a(0)>0.$
\end{itemize}
By rescaling, we can assume that $r=1$ and $V_{a}(x)>0$ for all $(x,a)\in [-1,1]\times A.$ In what follows we denote by~$\vartheta_{\varepsilon}$ the largest real root of $x(x^\mu-\varepsilon)=0$, i.e., 
\begin{equation}\label{sigma-equi}
\vartheta_{\varepsilon}=\left\{\begin{array}{ll} 0, &\text{if }\varepsilon\leqslant 0,\\[5pt]
\varepsilon^{1/\mu},&\text{if } \varepsilon\geqslant 0.
\end{array}\right.
\end{equation}
Observe then, see \refc{X}, that $(x,y)=(\vartheta_\varepsilon,0)$ is a hyperbolic saddle of $X$ for $\varepsilon\neq 0.$
We are interested in the Dulac time $\T(\,\cdot\,;\varepsilon,a)$  of the saddle-node unfolding~\refc{X} between the transverse sections $\{ y=1\}$ and $\{ x=1\}$. More concretely, see \figc{temps_SN}, for each $s>0$ small enough, we define $\T(s;\varepsilon,a)$ to be the time that spends the trajectory of $X$ starting at $(s+\vartheta_{\varepsilon},1)$ to arrive to $\{x=1\}$.

\begin{bigtheo}\label{A}
The Dulac time $\T(s;\varepsilon,a)$ of the saddle-node unfolding~\refc{X} between the transverse sections $\{ y=1\}$ and $\{ x=1\}$ verifies that
\[
 \lim_{(s, \varepsilon) \to (0^+,0)}\partial_s \T(s;\varepsilon, a)=-\infty
\] 
uniformly $($with respect to $a)$ on every compact subset of $A$.
\end{bigtheo}

In the next result we consider the family of dehomogenized Loud's centers \refc{loud0}, whose study is the main motivation of the present paper.

\begin{bigtheo}\label{regular}
Setting $a=(D,F)$, let $\{L_{a},a\in\R^2\}$ be the family of vector fields in~\refc{loud0} and consider the period function of the center at the origin. Then every $a=(D,F)\in (-1,0)\times\{0\}$ is a local regular value of the period function at the outer boundary of the period annulus. 
\end{bigtheo}

For a precise definition of local regular value we refer the reader to \cite[Definition 2.4]{MMV2}, but roughly speaking it means that no critical period bifurcates from these parameter values.

\section{Proofs of Theorems \ref{A} and \ref{regular}}

\begin{prooftext}{Proof of \teoc{A}.}
Let $y=y(x;s)$ be the trajectory of the vector field $x(x^\mu-\varepsilon)\partial_x-V_a(x)y\partial_y$ with initial condition $y(s+\vartheta_{\varepsilon};s)=1.$ Then there exist $s_0,\varepsilon_0>0$ small enough such that the Dulac time of \refc{X} between $\{y=1\}$ and $\{x=1\}$ is given by
\[
 \T(s;\varepsilon, a)=\int_{s+ \vartheta_{\varepsilon}}^1\left. \frac{U_a(x,y)}{x^\mu-\varepsilon}\right|_{y=y(x;s)}d\,x 
\] 
for all $s\in (0,s_0]$ and $\varepsilon\in [-\varepsilon_0,\varepsilon_0]$. Next, by applying the Weierstrass Division Theorem (see for instance \cite[Theorem 1.8]{Greuel} or \cite[Theorem 6.1.3]{Krantz}), we write 
\[
xU_a(x,y)=xU_a(x,0)+y\hat U_a(x,y),
\]
where $\hat U_a(x,y)$ is an analytic function on $(x,y,a)\in [-1,1]^2\times A$. Accordingly it turns out that $\T(s;\varepsilon,a)=\T_0(s;\varepsilon,a)+\T_1(s;\varepsilon,a)$ with
\[
 \T_0(s;\varepsilon,a)\!:=\int_{s+\vartheta_\varepsilon}^1\frac{U_a(x,0)}{x^\mu-\varepsilon}dx
 \text{ and }
 \T_1(s;\varepsilon,a)\!:=\int_{s+\vartheta_\varepsilon}^1\frac{y\hat U_a(x,y)}{x(x^\mu-\varepsilon)}dx.
\]
Note, and this is the key point, that $\T_1(s;\varepsilon,a)$ is the Dulac time of the saddle-node unfolding
\begin{equation*}
\hat X= \frac{1}{y\hat U_a(x,y)}\Big(x(x^\mu-\varepsilon)\partial_x-V_a(x)y\partial_y\Big),
\end{equation*}
which is in the hypothesis of \cite[Corollary B]{MMSV}. Thus, by applying that result with $\ell=k=1$, we obtain functions $c_{0}(\varepsilon,a)$ and $c_{1}(\varepsilon,a)$, satisfying that for every compact set $K_{a}\subset A$, there exists $\varepsilon_{1}>0$ such that $c_{0}$ and $c_{1}$ are continuous on 
$[-\varepsilon_{1},\varepsilon_{1}]\times K_{a}$ and 
\begin{equation}\label{eq1}
\T_1(s;\varepsilon,a)=c_{0}(\varepsilon,a)+c_{1}(\varepsilon,a)s+sh(s;\varepsilon,a),
\end{equation}
where the function $h$ in the remainder verifies $\lim_{s\to 0^+}h(s;\varepsilon,a)=0$ and $\lim_{s\to 0^+}s\partial _sh(s;\varepsilon,a)=0$ 
uniformly on $[-\varepsilon_{1},\varepsilon_{1}]\times K_{a}$. Observe on the other hand that 
\[
\partial_s \T_0(s;\varepsilon, a)=-\frac{U_a(s+ \vartheta_{\varepsilon},0)}{(s+ \vartheta_{\varepsilon})^\mu-\varepsilon}.
\]
Thus, the hypothesis $U_a(0,0)>0$ for all $a\in A$ and the fact that $\vartheta_\varepsilon$ tends to $0$ as $\varepsilon\to 0$, imply that for each compact set $K\subset A$ there exist positive constants $M,$ $s_0$ and $\varepsilon_0$ such that  $U_a(s+ \vartheta_{\varepsilon}, 0) > M$ for all $s \in (0,s_0] $, $\varepsilon \in [-\varepsilon_{0}, \varepsilon_{0}]$ and $a \in K$. Accordingly
\[
 \lim_{(s, \, \varepsilon) \to (0^+, \, 0)}\partial_s \T_0(s;\varepsilon, a)= -\infty
\] 
uniformly on every compact subset of $A.$ Consequently, taking \refc{eq1} also into account,
\[
 \partial_s\T(s;\varepsilon_a)=\partial_s \T_0(s;\varepsilon, a)+c_1(\varepsilon,a)+h(s;\varepsilon,a)+s\partial _sh(s;\varepsilon,a)\to -\infty
\]
as $(s,\varepsilon)\to(0^+,0)$ uniformly on $K_a.$ This concludes the proof of the result.
\end{prooftext}
\color{black}
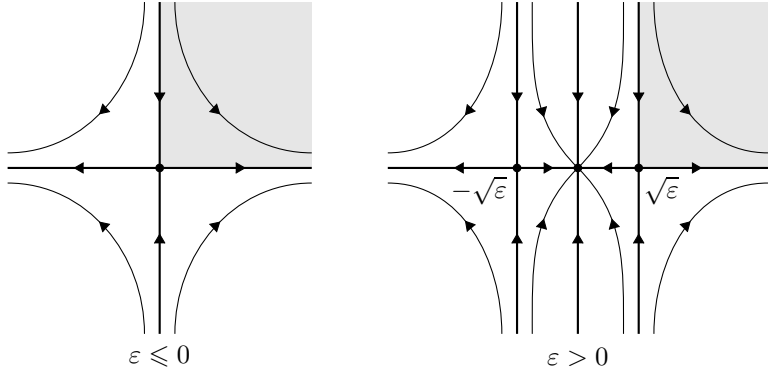
\begin{figure}[t]
\begin{center}
\begin{tikzpicture}
\path [fill=gray!20!white] (4.3,2.2) to (4.3,0) to (6,0) to (6,2.2) to (4.3,2.2);
\path [fill=gray!20!white] (-2,2.2) to (-2,0) to (0,0) to (0,2.2) to (-2,2.2);
\draw [fill]  (-2,0) circle [radius=0.05];
\draw [fill]  (3.5,0) circle [radius=0.05];
\draw [fill]  (2.7,0) circle [radius=0.05];
\draw [fill]  (4.3,0) circle [radius=0.05];
\draw [thick] (-4,0) to (0,0);
\draw [thick] (-2,-2.2) to (-2,2.2);
\draw (-1.8,2.2) to [out=-90,in=180] (0,.2);
\draw (-1.8,-2.2) to [out=90,in=180] (0,-.2);
\draw (-2.2,2.2) to [out=-90,in=0] (-4,.2);
\draw (-2.2,-2.2) to [out=90,in=0] (-4,-.2);
\draw [thick] (1,0) to (6,0);
\draw [thick] (2.7,-2.2) to (2.7,2.2);
\draw [thick] (4.3,-2.2) to (4.3,2.2);
\draw [thick] (3.5,-2.2) to (3.5,2.2);
\draw (2.9,2.2) to [out=-90, in=135] (3.5,0) to [out=-45,in=90] (4.1,-2.2);
\draw (4.1,2.2) to [out=-90, in=45] (3.5,0) to [out=-135,in=90] (2.9,-2.2);
\draw (4.5,2.2) to [out=-90,in=180] (6,.2);
\draw (4.5,-2.2) to [out=90,in=180] (6,-.2);
\draw (2.5,2.2) to [out=-90,in=0] (1,.2);
\draw (2.5,-2.2) to [out=90,in=0] (1,-.2);
\node at (-2,-2.5) {$\varepsilon\leqslant 0$};
\node at (3.5,-2.5) {$\varepsilon>0$};
\node at (4.6,-0.3) {$\sqrt{\varepsilon}$};
\node at (2.2,-0.3) {$-\sqrt{\varepsilon}$};
\begin{scope}[shift={(-1.3,0.8)},rotate=-42.5]
\draw[fill]  (0,.07) to (.12,0) to (0,-0.07);
\end{scope}
\begin{scope}[shift={(-1.3,-0.8)},rotate=42.5]
\draw[fill]  (0,.07) to (.12,0) to (0,-0.07);
\end{scope}
\begin{scope}[shift={(-2.7,0.8)},rotate=222.5]
\draw[fill]  (0,.07) to (.12,0) to (0,-0.07);
\end{scope}
\begin{scope}[shift={(-2.7,-0.8)},rotate=-222.5]
\draw[fill]  (0,.07) to (.12,0) to (0,-0.07);
\end{scope}
\begin{scope}[shift={(-1,0)},rotate=0]
\draw[fill]  (0,.07) to (.12,0) to (0,-0.07);
\end{scope}
\begin{scope}[shift={(-3,0)},rotate=180]
\draw[fill]  (0,.07) to (.12,0) to (0,-0.07);
\end{scope}
\begin{scope}[shift={(-2,1)},rotate=-90]
\draw[fill]  (0,.07) to (.12,0) to (0,-0.07);
\end{scope}
\begin{scope}[shift={(-2,-1)},rotate=90]
\draw[fill]  (0,.07) to (.12,0) to (0,-0.07);
\end{scope}
\begin{scope}[shift={(3.5,1)},rotate=-90]
\draw[fill]  (0,.07) to (.12,0) to (0,-0.07);
\end{scope}
\begin{scope}[shift={(3.5,-1)},rotate=90]
\draw[fill]  (0,.07) to (.12,0) to (0,-0.07);
\end{scope}
\begin{scope}[shift={(3.95,0)},rotate=180]
\draw[fill]  (0,.07) to (.12,0) to (0,-0.07);
\end{scope}
\begin{scope}[shift={(3.05,0)},rotate=0]
\draw[fill]  (0,.07) to (.12,0) to (0,-0.07);
\end{scope}
\begin{scope}[shift={(5,0)},rotate=0]
\draw[fill]  (0,.07) to (.12,0) to (0,-0.07);
\end{scope}
\begin{scope}[shift={(2,0)},rotate=180]
\draw[fill]  (0,.07) to (.12,0) to (0,-0.07);
\end{scope}
\begin{scope}[shift={(2.7,1)},rotate=-90]
\draw[fill]  (0,.07) to (.12,0) to (0,-0.07);
\end{scope}
\begin{scope}[shift={(2.7,-1)},rotate=90]
\draw[fill]  (0,.07) to (.12,0) to (0,-0.07);
\end{scope}
\begin{scope}[shift={(4.3,1)},rotate=-90]
\draw[fill]  (0,.07) to (.12,0) to (0,-0.07);
\end{scope}
\begin{scope}[shift={(4.3,-1)},rotate=90]
\draw[fill]  (0,.07) to (.12,0) to (0,-0.07);
\end{scope}
\begin{scope}[shift={(4.89,0.8)},rotate=-45]
\draw[fill]  (0,.07) to (.12,0) to (0,-0.07);
\end{scope}
\begin{scope}[shift={(4.89,-0.8)},rotate=45]
\draw[fill]  (0,.07) to (.12,0) to (0,-0.07);
\end{scope}
\begin{scope}[shift={(3.02,-0.8)},rotate=75]
\draw[fill]  (0,.07) to (.12,0) to (0,-0.07);
\end{scope}
\begin{scope}[shift={(3.02,0.8)},rotate=-75]
\draw[fill]  (0,.07) to (.12,0) to (0,-0.07);
\end{scope}
\begin{scope}[shift={(3.98,-0.8)},rotate=105]
\draw[fill]  (0,.07) to (.12,0) to (0,-0.07);
\end{scope}
\begin{scope}[shift={(3.98,0.8)},rotate=-105]
\draw[fill]  (0,.07) to (.12,0) to (0,-0.07);
\end{scope}
\begin{scope}[shift={(2.1,-0.8)},rotate=130]
\draw[fill]  (0,.07) to (.12,0) to (0,-0.07);
\end{scope}
\begin{scope}[shift={(2.1,0.8)},rotate=-130]
\draw[fill]  (0,.07) to (.12,0) to (0,-0.07);
\end{scope}
\end{tikzpicture}
\end{center}
\caption{Phase portrait of the orbital normal form $x(x^2-\varepsilon)\partial_x-V_a(x)y\partial_y$. 
The normalizing change of coordinates maps the period annulus of the Loud center to the gray quadrant.}\label{fig2}
\end{figure}
\begin{prooftext}{\bf Proof of Theorem \ref{regular}. } 
For the sake of convenience we reverse time in the original dehomogeneized Loud family \refc{loud0} and consider the vector field $-L_a$ instead. To study the saddle-node bifurcation that occurs at infinity we work in the projective plane $\RP^2$ and perform the change of coordinates
\[
 (z,w)=\left(\frac{1}{v},\frac{1-u}{v}\right).
\] 
The meromorphic extension of $-L_a$ in these coordinates is given by
\begin{align}\label{bXa}
\bar X_a\!:=\frac{1}{z}\Big(z\big(F+(D+1)z^{2}-&(2D+1)zw+Dw^{2}\big)\partial_z
\\\notag
&+w\big(-1+F+(D+1)z^{2}-(2D+1)zw+Dw^{2}\big)\partial_w\Big).
\end{align}
Some long but easy computations show that the local analytic change of coordinates given by
\begin{align}\label{cc}
&(x,y)=\Psi(z,w)\!:=\left(\frac{z}{\sqrt{g(z,w)}},\frac{w}{\sqrt{g(z,w)}}\right),\\
\intertext{where}\notag
&g(z,w)\!:=\frac{D}{2(F-1)(D+1)}w^2-\frac{(2D+1)}{(2F-1)(D+1)}wz+\frac{1}{2(D+1)},
\end{align}
brings the vector field $\bar X_a$ in \refc{bXa} to
\begin{align}\label{Xa}
&X_a=\frac{1}{xU_a(x,y)}\Big(x(x^2+2F)\partial_x+y(x^2+2F-2)\partial_y\Big),\\
\intertext{where}\notag
&U_a(x,y)\!:=\left(\frac{(2D+1)}{2(2F-1)}xy-\frac{D}{4(F-1)}y^2+\frac{D+1}{2}\right)^{-1/2}.
\end{align}
Indeed, one can verify
%
%
that $\Psi^*X_a=(D\Psi)^{-1}(X_a\circ\Psi)=\bar X_a.$ For reader's convenience let us briefly explain how we obtain this normalizing change of coordinates. The idea arises from the fact that $\bar I(z,w)=\frac{w}{z}\left( 1+2F\frac{g(z,w)}{z^2}\right)^{-\frac{1}{2F}}$ is a first integral of $\bar X_a$ for $F\neq 0$, and that the change of coordinates in \refc{cc} brings it to $I(x,y)=\frac{y}{x}\left(1+\frac{2F}{x^2}\right)^{-\frac{1}{2F}}$. Thus, since the $1$-form $dI$ is proportional to $x(x^2+2F)dy-y(x^2+2F-2)dx$, we deduce that the coordinate change brings $\bar X_a$ to 
\[
 K_a(x,y) \Big(x(x^2+2F)\partial_x+y(x^2+2F-2)\partial_y\Big),
\]
so that the problem reduces to find this factor $K_a.$

We are now in position to apply \teoc{A} because, taking
\[
 \mu\!:=2,\;\varepsilon\!:=-2F\text{ and }V_a(x)=2-2F-x^2,
\]
observe that we can write the vector field in \refc{Xa} as
\[
X_a=\frac{1}{xU_a(x,y)}\Big(x(x^\mu-\varepsilon)\partial_x-yV_a(x)\partial_y\Big).
\]
Note moreover that in these normalizing coordinates the period annulus of the center at the origin of $-L_a$ is contained in the quadrant $\{x \geqslant \vartheta_{\varepsilon},y\geqslant 0\},$ see \figc{fig2}, where $\vartheta_{\varepsilon}$ is given in~\refc{sigma-equi} with $\mu=2$. Following the notation in \teoc{A} we also take 
\[ 
 A\!:=(-1,0)\times(-1/2,1/2).
\]  
Then $U_a(0,0)>0$ and $V_a(0)>0$ for all $a\in A$. Furthermore, 
working on any compact subset $K_{a}$ of~$A$, we see that the Taylor series of $U_{a}(x,y)$ at $(0,0)$ is absolutely convergent for $(x,y)\in [-r,r]^2$ for some $r>0$ depending only on $K_a.$ By rescaling the local coordinates $(x,y)$ we can assume that $r=1.$ This will change $U_a$, $V_a$ and $\varepsilon$ in terms of $a$ but it will be clear that the proof does not depend on their particular expression, provided that the new $\varepsilon$ tends to zero as $F\to 0,$ which one can verify that this is the case. 

Next we proceed to study the period function of the center near the polycycle at the boundary of its period annulus. To this end we first note that the vector field in~\refc{loud0} is invariant with respect to the symmetry $(u,v)\to (u,-v)$, and so is $-L_a.$ Consequently, see \figc{disc}, its period function is twice the time that the solutions of $-L_a$ spend for going from $\Sigma_{1}\!:=\{u \approx -\infty , v=0\}$ to $\Sigma_{2}\!:=\{u \approx 1, v=0\}$.
\begin{figure}[t]
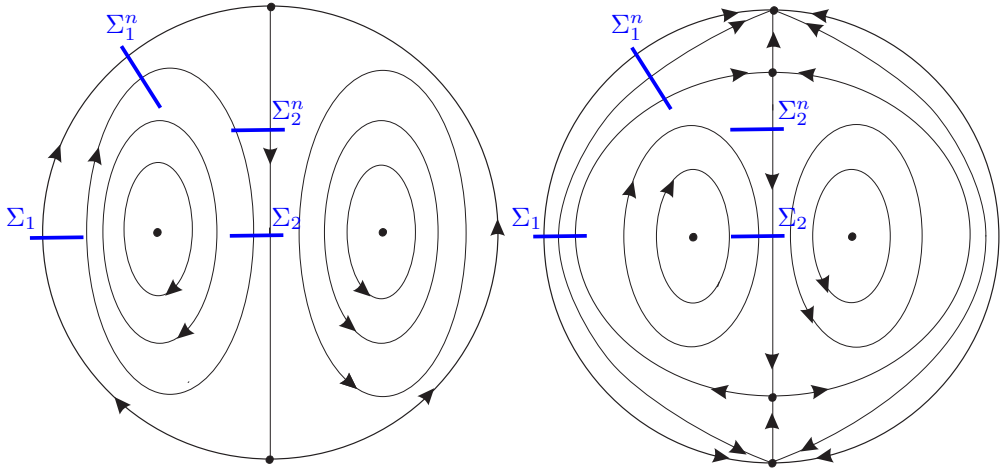

 \centering
 \begin{lpic}[l(0mm),r(0mm),t(0mm),b(5mm)]{fig4(1.75)}
   \lbl[l]{-1.5,19;$\blue{\Sigma_1}$}   
   \lbl[l]{36.5,19;$\blue{\Sigma_1}$} 
   \lbl[l]{18.5,19;$\blue{\Sigma_2}$}   
   \lbl[l]{56.5,19;$\blue{\Sigma_2}$} 
   \lbl[l]{18.5,27;$\blue{\Sigma_2^n}$}   
   \lbl[l]{56.5,27;$\blue{\Sigma_2^n}$}  
   \lbl[l]{6,33.5;$\blue{\Sigma_1^n}$}   
   \lbl[l]{44,33.5;$\blue{\Sigma_1^n}$} 
 \end{lpic}
\caption{Phase portrait of the vector field $-L_a$ for $D\in (0,1)$ and $F\approx 0$ in the Poincaré disc. On the left, for $F\geqslant 0$, and on the right, for $F<0$. By symmetry, the period of a periodic orbit is twice the time that spends the solution for going from $\Sigma_1$ to $\Sigma_2$. The auxiliary transverse sections $\Sigma_1^n$ and $\Sigma_2^n$ are defined by means of the normalizing change of coordinates $\Psi$ that brings, locally, the saddle-node unfolding at infinity to the normal form in \refc{X}, so that \teoc{A} applies.}\label{disc}
\end{figure}
In order to study this, let us say, half period function we introduce two auxiliary transverse sections near the saddle-node bifurcation, $\Sigma_{1}^{n}\!:=\Psi^{-1}(\{y=1\})$ and $\Sigma_{2}^{n}\!:=\Psi^{-1}(\{x=1\})$, 
parameterized by $s\mapsto \Psi^{-1}(s+\vartheta_{\varepsilon},1)$ and $s\mapsto \Psi^{-1}(1,s),$ respectively. 
Here~$\Psi$ is the (local) normalizing change of coordinates given in \refc{cc} and we work with the projective coordinates $(z,w)$. Then we define $T(s;a)$ to be half of the period of the periodic orbit $\gamma_{a,s}$ of $-L_a$ passing through the  point $\Psi^{-1}(s+\vartheta_{\varepsilon},1)\in\Sigma_1^n$ and we decompose it as 
\[
 T(s;a)=T_{1}(s;\varepsilon,a)+\T(s;\varepsilon,a)+T_{2}\big(\mathscr D(s;\varepsilon,a);\varepsilon,a\big),
\]
where (see \figc{disc} again)
\begin{itemize}
\item $T_{1}(s;\varepsilon,a)$ is the time that spends $\gamma_{a,s}$ for going from $\Sigma_1$ to $\Sigma_1^n$,
\item $\T(\,\cdot\,;\varepsilon,a)$ and $\mathscr D(\,\cdot\,;\varepsilon,a)$ are the Dulac time and Dulac map from $\Sigma_1^n$ to $\Sigma_2^n$, respectively, 
\item and $T_2(\,\cdot\,;\varepsilon,a)$ is the transition time form $\Sigma_2^n$ to $\Sigma_2$.
\end{itemize}
(Here the dependence on $\varepsilon$ is redundant because it is a function of $a=(D,F)$ but we keep it to be consistent with the notation of \teoc{A}.) We next study each one of these summands. To this end, given any compact subset $K_a$ of~$A$, we denote by $\mathcal I(K_a)$ the space of functions $h(s;a)$, analytic on $s\in (0,s_{0}),$ verifying  
\[
 \lim_{s\to 0^+}h(s;a)=0\text{ and }\lim_{s\to 0^+}s\partial_{s}h(s;a)=0\text{ uniformly on $K_a.$}
\]
It is clear that $\mathcal I(K_a)$ is stable under addition and multiplication. 

Let us observe first that we can write $T_1(s;a)=f(s+\vartheta_{\varepsilon},a)$ where $f$ is an analytic function at $\{0\}\times A$, whereas the transition time $T_2(s;a)$ is analytic at 
$\{0\}\times A.$ Accordingly, for $i=1,2$, we have that 
\[
 T_{i}(s;a)=c_{i,0}(a)+c_{i,1}(a)s+sh_{i}(s;a)\text{ with $h_{i}\in\mathcal I(K_a)$}
\]
and where $c_{i,0}$ and $c_{i,1}$ are continuous functions on $K_a$. On the other hand, since $y=\mathscr D(x;\varepsilon,a)$ is a trajectory of the vector field $x(x^2-\varepsilon)\partial_x+V_a(x)y\partial_y$,
by applying assertion $(b)$ in \cite[Corollary~A]{MMSV} with $\{\mu=2,\ell=k=1,\lambda=2-2F)\}$, we deduce that $\mathscr D(s;\varepsilon,a)=s h_{0}(s;a)$ with $h_{0}\in\mathcal I(K_a).$ Consequently $T_{2}\big(\mathscr D(s;\varepsilon,a);\varepsilon,a\big)=c_{2,0}(a)+s\hat h_2(s;a)$ where $\hat h_{2}(s)\!:=c_{2,1}h_{0}(s)+h_{0}(s)h_{2}\bigl(sh_{0}(s)\bigr),$ and it can be easily checked that $\hat h_{2}\in\mathcal I(K_a).$ Summing up we can write 
\[
 T(s;a)=\T(s;\varepsilon,a)+c_0(a)+c_1(a)s+sh(s;a)
\]
where $c_0\!:=c_{1,0}+c_{2,0}$ and $c_1\!:=c_{1,1}$ are continuous functions on $K_a$ and $h\!:=h_1+h_2+\hat h_2\in\mathcal I(K_a)$. Taking this into account, since $a=(D,F)$, the application of \teoc{A} to $\T(s;\varepsilon,a)$ shows that the derivative $\partial_s T(s;D,F)$ tends to $-\infty$ as $(s,D,F)\to (0^+,D_0,0)$ for every $D_0\in(-1,0)$. Consequently there exists $\delta>0$ such that $\partial_s T(s;D,F)\neq 0$ for all $(s,D,F)$ with $s\in (0,\delta),$ $|D-D_0|<\delta$ and $|F|<\delta.$ Since the period of the periodic orbit $\gamma_{s,a}$ is $2T(s;D,F)$, this concludes the proof of the result.
\end{prooftext}

\bibliographystyle{plain}

\begin{thebibliography}{99}

\bibitem{Chicone1} C. Chicone, {\it The monotonicity of the period function for planar Hamiltonian vector fields}, J. Differential Equations \textbf{69} (1987) 310--321.

\bibitem{Chicone2} C. Chicone and M. Jacobs, {\it Bifurcation of critical periods for plane vector fields}, Trans. Amer. Math. Soc. \textbf{312} (1989) 433--486.

\bibitem{Chicone3} C. Chicone, review in MathSciNet, ref. 94h:58072.

\bibitem{Colin} C. Christopher and J. Devlin, {\it On the classification of Liénard systems with amplitude-independent periods}, J. Differential Equations \textbf{200} (2004) 1--17.

\bibitem{CMV} A. Cima, F. Mañosas and J. Villadelprat, {\it Isochronicity for several classes of Hamiltonian systems}, J.~Differential Equations \textbf{157}  (1999) 373--413.

\bibitem{CG}  W.A. Coppel and L. Gavrilov, {\it The period function of a Hamiltonian quadratic system,} Differential Integral Equations {\bf 6} (1993) 1357--1365.

\bibitem{ADL} F. Dumortier, J. Llibre and J.C. Artés, ``Qualitative theory of planar differential systems'',
Universitext, Springer-Verlag, Berlin, 2006.

\bibitem{Fonda} A. Fonda, M. Sabatini and F. Zanolin, {\it Periodic solutions of perturbed Hamiltonian systems in the plane by the use of the Poincaré-Birkhoff Theorem}, Topol. Methods Nonlinear Anal. \textbf{40} (2012) 29--52.

\bibitem{Pugh} J.-P. Françoise and C. Pugh, {\it Keeping track of limit cycles},
J. Differential Equations \textbf{65} (1986) 139--157. 

\bibitem{GGV} A. Gasull, A. Guillamon and J. Villadelprat, {\it The period function for second-order
quadratic ODEs is monotone,} Qual. Theory Dyn. Syst. \textbf{4} (2004) 329--352.

\bibitem{GMM} A. Gasull, V. Mañosa and F. Mañosas, {\it Stability of certain planar unbounded polycycles}, J. Math. Anal. Appl. \textbf{269} (2002) 332--351.

\bibitem{Gasull} A. Gasull, C. Liu and J. Yang, {\it On the number of critical periods for planar polynomial systems of arbitrary degree}, J. Differential Equations \textbf{249} (2010) 684--692.

\bibitem{GV10} M. Grau and J. Villadelprat, {\it Bifurcation of critical periods from Pleshkan's isochrones,}
J. London Math. Soc. \textbf{81} (2010) 142--160. 

\bibitem{Greuel} G.-M. Greuel, C. Lossen and E. Shustin, ``Introduction to singularities and deformations'', Springer Monogr. Math., Springer, Berlin, 2007.

\bibitem{GH} J. Guckenheimer and P. Holmes, ``Nonlinear oscillations,
dynamical systems, and bifurcations of vector fields'', Appl. Math. Sci. \textbf{42}, Springer-Verlag, New York, 1983.

\bibitem{Krantz} S. G. Krantz and H. R. Parks, ``A Primer of Real Analytic Functions'', Birkh\"auser Advanced Texts, Birkh\"auser Basel, 2002.

\bibitem{MMV} P. Marde\v si\'c, D. Mar{\'\i}n and J. Villadelprat, {\it On the time function of the Dulac map for families of meromorphic vector fields,} Nonlinearity \textbf{16} (2003),  855--881.

\bibitem{MMV2} P. Marde\v si\'c, D. Mar{\'\i}n and J. Villadelprat, {\it The period function of reversible quadratic
centers,} J. Differential Equations \textbf{224} (2006), 120--171.

\bibitem{MMV3} P. Marde\v si\'c, D. Mar{\'\i}n and J. Villadelprat, {\it Unfolding of resonant saddles and the Dulac time,} Discrete and Continuous Dynamical Systems, \textbf{21} (2008), 1221--1244.

\bibitem{MMSV}  P. Marde\v si\'c, D. Mar{\'\i}n, M. Saavedra and J. Villadelprat, \emph{Unfoldings of saddle-nodes and their Dulac time}. J. Differential Equations \textbf{261} (2016), 6411 - 6436.

\bibitem{MV} D. Mar\'{\i}n and J. Villadelprat, {\it On the return time function around monodromic polycycles,} J. Differential Equations \textbf{228} (2006) 226--258.

\bibitem{Roussarie} R. Roussarie, ``Bifurcation of planar vector fields and Hilbert's sixteenth problem'', Progr. Math. \textbf{164}, Birkhäuser Verlag, Basel, 1998.

\bibitem{Saa20} M. Saavedra, {\it Dulac time of a resonant saddle in the Loud family} J. Differential Equations {\bf 269} (2020) 7705--7729.

\bibitem{V07} J. Villadelprat, {\it On the reversible quadratic centers with monotonic period function}, Proc. Amer. Math. Soc. \textbf{135} (2007) 2555--2565.

\bibitem{Xingwu} C. Xingwu and W. Zhang, {\it Isochronicity of centers in a switching Bautin system},
J. Differential Equations \textbf{252} (2012) 2877--2899. 

\bibitem{Ye} Yan-Qian Ye {\it et al}, ``Theory of limit cycles'', Transl. Math. Monogr. \textbf{66}, American Mathematical Society, Providence, RI, 1986.

\bibitem{Zhao} Y. Zhao, {\it On the monotonicity of the period function of a quadratic system},
Discrete Contin. Dyn. Syst. \textbf{13} (2005) 795--810.



\end{thebibliography}

\end{document}